\documentclass[12pt]{article}

\usepackage{graphicx}
\usepackage{amsmath,amssymb}

\newcommand{\E}{\mathbb E}
\newcommand{\N}{\mathbb N}
\newcommand{\R}{\mathbb R}

\newcommand{\dd}{\text{\rm d}}
\newcommand{\x}{\mathbf x}
\newtheorem{theorem}{Theorem}
\newtheorem{lemma}{Lemma}

\newtheorem{proposition}{Proposition}
\newtheorem{remark}{Remark}

\begin{document}



\title{Approaches to asymptotics for $U$-statistics of Gibbs facet processes}


\author{Jakub Ve\v{c}e\v{r}a, Viktor Bene\v{s}\\
Charles University in Prague, Faculty of Mathematics and Physics,\\ Department of Probability and Mathematical Statistics,\\ Sokolovsk\'{a} 83, 18675 Praha 8, Czech Republic}
\maketitle
\subsection*{Abstract}
It is shown how the central limit theorem for $U$-statistics of spatial Poisson point processes can help to derive the central limit theorem for $U$-statistics of a Gibbs facet process from stochastic geometry. A full-dimensional submodel enables a simpler approach to its investigation. Finally the general situation is studied and the asymptotics with increasing intensity is described.

\noindent{\bf Keywords:} central limit theorem; Gibbs facet process; $U$-statistics

\noindent{\bf MSC:} {60D05}

\section{Introduction}

Central limit theorems for $U$-statistics of spatial Poisson processes with increasing intensity were derived based on Malliavin-Stein method in \cite{RefR}. Later in \cite{RefLT} a multivariate version of this result was achieved using the Malliavin-Stein method and alternatively the moment and/or the cumulant method. The aims to extend developments of this type to functionals of a wider class of spatial processes, e.g. Gibbs processes \cite{RefCh}, were initiated by \cite{RefX}.

In \cite{RefV} we introduce facet processes in arbitrary Euclidean dimension. These are finite Gibbs type processes of compact subsets of hyperplanes in which repulsive interactions enter by penalizing intersections of facets. In dimension 2 or 3, cf. \cite{RefB}, facet process is the segment, surface process, respectively, which may serve for modeling real data from biology or materials research. The exponential type density w.r.t. Poisson process has been used earlier in a planar disc model by \cite{RefK}, \cite{RefM}. For both mathematical and statistical purposes limit behavior of functionals of the process is of interest. In \cite{RefV} natural $U$-statistics of the model are studied in case when the intensity of the reference Poisson process tends to infinity. In \cite{RefW} the multivariate central limit theorem for this system of $U$-statistics of the facet process is derived.

In the present paper two results concerning the facet process are derived which complement and improve \cite{RefV} and \cite{RefW}. First we show how the full-dimensional submodel can be investigated in a simpler way which enables to track the distribution on the way to the limit case. Secondly the proof of the multivariate central limit theorem for the vector of $U$-statistics of the Gibbs facet process is naturally shortened using the recent result for the Poisson process from \cite{RefLT}. 
 
\section{The facet process}
We call facets compact subsets of hyperplanes in $\R^d$ with a given shape. In the space $Y=B\times (0,b]\times {\mathbb S}^{d-1}$ a compact window $B$ is the set of facet centres, $(0,b]$ is an interval of limited sizes of facets, ${\mathbb S}^{d-1}$ is the hemisphere of normal orientations. Since the investigation of a general Gibbs type model for randomly dispersed facets with interactions is hardly tractable, recently in \cite{RefV} and also in the present paper we have been studying a special case
\begin{equation}\label{spc}Y=[0,b]^{d} \times \{2b \} \times \{ e_{1},\ldots,e_{d} \},\end{equation} i.e. the size is fixed and $d$ orientations correspond to the basic orthonormal system of vectors in $\R^d.$ $Y$ is isomorphic with the space of facets, i.e. for a point $((z_{1},\dots,z_{d}),2b,e_{l})\in Y,\; l\in [d]=\{1,\dots ,d\}$ a facet is defined to be the set
$$ \{(x_{1},\ldots,x_{d}) \in \mathbb{R}^{d}, x_{l}=z_{l};|x_{i}-z_{i}|\leq b,i \in [d] \setminus \{l \} \}.$$ In this setting all non-parallel facets intersect and denoting $\mathbb{H}^k$ the Hausdorff measure of order $k$ in $\R^d,$ we have bounds for the measure of intersection of $c$ mutually non-parallel facets \begin{equation}\label{bnds}b^{d-c} \leq \mathbb{H}^{d-c}(\cap_{i=1}^{c}y_{i}) \leq (2b)^{d-c} .\end{equation}
Consider a finite measure $\lambda$ on $Y$ of a form
\begin{equation}\label{ori}
\lambda (\dd x)=\lambda (\dd (z,r,\phi)) = \chi (z) \dd z \delta_{2b}(\dd r) \frac{1}{d} \sum_{i=1}^{d} \delta_{e_{i}} (\dd \phi),
\end{equation}
where $\delta $ is the Dirac measure and $\chi $ is a bounded intensity function of facet centres.
Further let  $({\mathbf N},{\mathcal N})$ be a measurable space of integer-valued finite measures on $Y$, where each atom has measure one. Here ${\mathcal N}$ is the smallest $\sigma$-algebra which makes the mappings $\x\mapsto \x(A)$ measurable for all Borel sets $A\subset Y$ and all $\x\in {\mathbf N}.$ 
An alternative meaning of $\x$ is being the support of the measure, a finite point set.
 
A $U$-statistic of order $k$ is a measurable function $F$ on ${\mathbf N}$ given by a formula \begin{equation}\label{ust}F(\mathbf{x})=\sum_{(x_{1},\ldots,x_{k}) \in \mathbf{x}^{k}_{\neq}}f(x_1,\dots ,x_k),\end{equation} where $f\in L^1(\lambda^k)$ is called a driving function. $f$ is symmetric which means invariant with respect to permutations of its variables. In (\ref{ust}) we sum over ordered $k$-tuples of distinct points from $\mathbf{x}.$

Let $(\Omega ,{\cal A}, P)$ be a probability space. We denote for $a\geq 1$$$\eta_{a}:\,(\Omega ,{\cal A}, P)\rightarrow ({\mathbf N},{\mathcal N})$$ a finite Poisson process of facets with intensity measure $a\lambda .$  We will consider a system of $U$-statistics of order $j$ \begin{equation}\label{gecka}G_j,\; j=1,\dots ,d,\end{equation} with driving functions $g^{(j)}$ defined as
\begin{align}
\nonumber g^{(j)}(x_{1},\ldots,x_{j})=\frac{1}{j!}\mathbb{H}^{d-j} (\cap_{i=1}^{j} x_{i} ),\; j=1,\dots ,d.
\end{align} 
Further from \cite{RefR} $$g^{(j)}_n(y_1,\dots ,y_n)=\binom{j}{n}\int g^{(j)}(y_1,\dots ,y_n,x_1,\dots ,x_{j-n})\lambda(\dd (x_1,\dots ,x_{j-n})),\; n\leq j,$$ using notation $\lambda(\dd (x_1,\dots ,x_{k}))=\lambda (\dd x_1)\dots\lambda (\dd x_k).$
We will standardize the vector (\ref{gecka}) of $U$-statistics in a form \begin{equation}\label{stdz}\tilde{G}_{j}(\eta_a) = \frac{ G_{j}(\eta_{a})- \E G_{j}(\eta_{a})}{a^{j-\frac{1}{2}}},\; 1 \leq j  \leq d.\end{equation} The following asymptotic results for Poisson processes when $a\rightarrow\infty $
follow from \cite{RefLT}. The asymptotic covariances are \begin{equation}\label{ascv}C_{ij}=\lim_{a\rightarrow\infty }cov(\tilde{G}_{i}(\eta_a),\tilde{G}_{j}(\eta_a))=\langle g_1^{(i)},\,g_1^{(j)}\rangle_1 ,\end{equation} where $\langle\cdot ,\cdot\rangle_p$ is the inner product in $L^2(\lambda^p).$
\begin{theorem}\label{lpst}Let $M_d$ be a $d$-dimensional centered Gaussian random vector with covariance matrix $(C_{ij}),\,i,j=1,\dots ,d.$ Then
$$(\tilde{G}_{1}(\eta_a),\dots ,\tilde{G}_{d}(\eta_a)) \rightarrow M_d$$ in distribution.
\end{theorem} The integrability assumptions of the central limit theorem as stated in Proposition 5.1 in \cite{RefLT} are fulfilled in our setting since the space $Y$ is bounded and the intensity measure is finite.

In the following we will study the facet process $\mu_{a}$ with a density \begin{equation}\label{dnst}p(\mathbf{x})=c_{a} \exp \left( \sum_{i=1}^{d} \nu_{i} G_{i}(\mathbf{x}) \right)\end{equation} with respect to $\eta_{a}$, where $a \geq 1$, $\nu_{i},\,i=1,\dots ,d,$ are real parameters and $ c_{a}$ is the normalizing constant. Fulfilling of condition $\nu_{i} \leq 0 ,i=2,\ldots,d$ assures that $p \in L^{1}(P_{\eta_a}) \cap L^{2}(P_{\eta_a}),$ where $P_{\eta_a}$ is the probability distribution of $\eta_a.$ We will use the notion of a submodel of order $l$ $$\mu_{a}^{(l)},\;1\leq l\leq d,$$ where in (\ref{dnst}) we have $\nu_{j}=0,\: j\neq l,\;\nu_l<0.$ The properties of submodels of the order higher than $1$ will be explored, since $\mu_{a}^{(1)}$ is a Poisson process.

\section{The full-dimensional submodel}
Among the $U$-statistics (\ref{gecka}) $G_d$ plays a special role and in this section we restrict our attention to it. Without much loss of generality we assume $\chi (z)\equiv 1.$ In $g^{(d)}$ we have $\mathbb{H}^0$ which is a counting measure and therefore $G_d$ depends only on orientations of facets and not on their locations. The submodel $\mu_a^{(d)}$ will be studied, writing $\nu_d=\nu .$ Consider a map 
\begin{equation}\label{sgm}\theta:\, ({\mathbf N},{\mathcal N})\rightarrow\mathbb{N}_{0}^d;\quad \x\mapsto (\theta_1(\x ),\dots ,\theta_d(\x )),\end{equation} where $\mathbb{N}_{0}=\mathbb{N}\cup\{0\}.$
 The $U$-statistic $G_d(\x ),\; \x\in {\mathbf N},$ can be expressed by means of variables $\theta_i (\x )$ which correspond to the numbers of facets in ${\mathbf x}$ which have orientation $e_i \; i=1,\dots ,d.$ We have that \begin{equation}\label{kack}G_d({\mathbf x})=\frac{1}{d!}\sum_{(x_{1},\ldots,x_{d}) \in \mathbf{x}^{d}_{\neq}}{\mathbb H}^0(\cap_{j=1} ^d x_j)=\prod_{i=1}^d\theta_i(\x ).\end{equation} A discrete probability distribution $$\pi= P_{\mu_a^{(d)}}\theta^{-1}$$ is defined on $\mathbb{N}_{0}^d,\;\pi (k_1,\dots ,k_d)$ is symmetric. 
\begin{proposition}For the distribution $\pi$ and $(k_1,\dots ,k_d)\in\mathbb{N}_{0}^d$ it holds \begin{equation}\label{targ}
\pi (k_1,\dots ,k_d)\propto \frac{A^{\sum k_i} }{k_1!\dots k_d!}\exp\left(\nu_d\prod_{i=1}^dk_i\right),\end{equation}
where $A=\frac{ab^d}{d}$ and the symbol $\propto$ means proportionality.
\end{proposition}
{\bf Proof:} We use the Radon-Nikodym theorem $$P_{\mu_a^{(d)}}(\dd \x )=p(\x )P_{\eta_a}(\dd\x ).$$ Random variables $\theta_i(\eta_a)$ are Poisson distributed with mean $$a\lambda (\{e_i\})=A,$$ where we used formula (\ref{ori}). Also here $\theta_i(\eta_a)$ are independent and thus $$P_{\eta_a}\theta^{-1}(k_1,\dots ,k_d)=e^{-dA}\frac{A^{\sum k_i} }{k_1!\dots k_d!}.$$ Using $\theta_i(\x)=k_i$ and $$p(\x )\propto\exp (\nu \prod_{i=1}^d\theta_i(\x ))$$ we obtain the result.
\hfill $\Box $

We will show that asymptotically for $a\rightarrow\infty $ the distribution $\pi $ tends to be concentrated on the set $\mathbb{N}_0^d\setminus\mathbb{N}^d.$
\begin{proposition} We have \begin{equation}\label{vta}\lim_{a\rightarrow\infty} \sum_{k_1=1}^\infty \dots\sum_{k_d=1}^\infty\pi(k_1,\dots ,k_d)=0.\end{equation}
\end{proposition}
{\bf Proof:} From Theorem 3 in \cite{RefV} $$\lim_{a\rightarrow\infty}\E G_d(\mu_a^{(d)})=0.$$ In our setting this means that \begin{equation}\label{pp}\E\left[\prod_{i=1}^d\theta_i(\mu_a^{(d)})\right]=\sum_{k_1=0}^\infty\dots\sum_{k_d=0}^\infty\left(\prod_{i=1}^dk_i\right)\pi (k_1,\dots ,k_d)\rightarrow 0,\;\;a\rightarrow\infty.\end{equation} The right hand side is zero if some $k_i=0$ and otherwise it is not smaller than the expression in (\ref{vta}). This implies the assertion of the Theorem.
\hfill $\Box $
 
\section{Multivariate central limit theorem}\label{sec}
In the following we study the asymptotic behavior in general situation for any submodel $\mu_a^{(c)}$ and a vector of $U$-statistics $G_j,\;2\leq c\leq d,\;1\leq j\leq d.$ We will focus on the central limit theorem since it has been already shown in \cite{RefV} that  

$\begin{array}{ll}\lim_{a\rightarrow\infty}\E G_j(\mu_a^{(c)})&=0,\; j\geq c,\\ &>0,\; j<c.\end{array}$

Let for an hereditary density $p,$ i.e. satisfying $p(\x )>0\Rightarrow p(\tilde{\x})>0$ whenever $\tilde{\x}\subset\x,$ $$\lambda^*_n(u_1,\dots ,u_n;\x)=\frac{p(\x\cup \{u_1,\dots ,u_n\})}{p(\x)},$$ $u_1,\dots ,u_n\in Y\setminus\{\x\}
$ distinct, be the conditional intensity of $n$-th order of $\mu_a ,\;\lambda^*_0\equiv 1.$ We observe that $\lambda^*_n$ is symmetric in the variables $u_1,\dots ,u_n.$ The expectation of conditional intensity \begin{equation}\label{coref}\rho_n(u_1,\dots ,u_n;\mu_a )={\mathbb E}\lambda^*(u_1,\dots ,u_n;\mu_a )={\mathbb E}[p(\eta_a\cup\{u_1,\dots ,u_n\})],\end{equation} is called $n-$th correlation function of the facet process $\mu_a,$ analogously for $\mu_a^{(c)}.$ 
\begin{lemma}\label{KorelFce}Let $p\in \mathbb{N},\; 2\leq c\leq d.$ Then there exist $R,S>0$ such that for any $a>0$ and $\{x_1,\dots ,x_p\}\subset Y$
we have \begin{align}
\nonumber \left\vert \rho_{p}(x_{1},\ldots,x_{p};\mu_{a}^{(c)}) - \frac{\binom{d-k}{d-c+1}}{\binom{d}{d-c+1}} \right\vert < Re^{-Sa},
\end{align} where $k$ is number of distinct facet orientations among $\{x_{1},\ldots,x_{p} \}.$
\end{lemma}
\noindent{\bf Sketch of the proof:}
First consider $ p \leq c$ and let facets $x_{1},\ldots,x_{p}$ have orientations $e_{1},\ldots,e_{p}$ (without loss of generality).
It holds \cite{RefB} $$\rho_{p}(x_{1},\ldots,x_{p};\mu^{(c)}_{a})=\frac{\E \exp(\nu_{c} G_{c}(\eta_{a}\cup \{x_{1},\ldots,x_{p}\} ))}{\E \exp(\nu_{c} G_{c}(\eta_{a} ))} =$$
\begin{align}
\label{zac}=\frac{\sum_{n=0}^{\infty}\frac{a^{n}}{n!} \int_{Y^{n}} \exp\left(    \nu_{c}G_{c} (\{ u_{1},\ldots,u_{n},x_{1},\ldots,x_{p} \})
 \right) \lambda^n (\dd (u_{1},\ldots,u_{n}))}{\sum_{n=0}^{\infty}\frac{a^{n}}{n!} \int_{Y^{n}} \exp\left(    \nu_{c}G_{c} (\{ u_{1},\ldots,u_{n} \}) \right) \lambda^n (\dd (u_{1},\ldots,u_{n}))}.
\end{align} Let $\mathbf{\theta }=(\theta_1,\dots ,\theta_d)=\theta(\{u_1,\dots ,u_n\}),$ cf. (\ref{sgm}). Define  \begin{equation}\label{rr}
A(c,q,s,\mathbf{\theta }) = \sum_{\substack{F \subset [s] \\
c-q \leq |F| \leq c \\
|F \cup [q]|\geq c
}} \prod_{j \in F } \theta_{j}.
\end{equation}$A(c,p,d,\mathbf{\theta })$ is the number of intersections of $c$-tuples of the facets among all $n+p$ facets with orientations of $u_1,\dots ,u_n$ described by $\mathbf{\theta }$ and orientations of $x_1,\dots ,x_p$ equal to $e_{1},\ldots,e_{p},$ respectively. For $\mathbf{n} =(n_1,\dots ,n_d)\in\mathbb{N}_{0}^d$ put
 \begin{equation}\label{bck}B(c,Q,a,q,s)=\sum_{n_{1}=0}^{\infty} \ldots \sum_{n_{d}=0}^{\infty} \frac{a^{n_{1} + \ldots +\ n_{d}}}{n_{1}! \ldots n_{d}!} \exp \left( \nu_{c} Q^{d-c} A(c,q,s,\mathbf{n}) - a(c-1)\right).\end{equation}
Substituting $a$ for $\frac{aT}{d},$ where $T= \int_{[0,b]^{d}} \chi(z) \dd z$
we obtain bounds for (\ref{zac}) by means of (\ref{bnds}) and (\ref{ori}):\begin{equation}\label{dbs}\frac{B(c,2b,a,p,d)}{B(c,b,a,0,d)}
\leq\rho_{p}(x_{1},\ldots,x_{p},\mu^{(c)}_{a}) \leq \frac{B(c,b,a,p,d)}{B(c,2b,a,0,d)}.
\end{equation} Using techniques analogous to those in \cite{RefV} it can be shown that for $q \leq c \leq s$, $c \geq 2$ and any $Q>0$ there exist $R,S>0$ such that
\begin{align}\label{cvg}
\left\vert B(c,Q,a,q,s) - \binom{s-q}{s-c+1} \right\vert <R e^{-Sa}. 
\end{align}Thus the limits of the bounds on both sides in (\ref{dbs}) are the same, also for other cases $(p>c).$ The rate of convergence (\ref{cvg}) extends to the fractions in (\ref{dbs}). 
\hfill $\Box $

The next step is the evaluation of moments of $U$-statistics.
We can use a short expression for moment formulas using diagrams and partitions, see \cite{RefP}, \cite{RefLT}. Let $\tilde{\prod}_k$ be the set of all  partitions $\{J_i\}$ of $[k]=\{1,\dots ,k\},$ where $J_i$ are disjoint blocks and $\cup J_i=[k].$ For
$k=k_1+\dots +k_m$ and blocks $$J_i=\{ j: k_1+\dots +k_{i-1}< j\leq k_1+\dots +k_i\},\; i=1,\dots ,m,$$ consider the partition $\pi =\{J_i,\;1\leq i\leq m\}$ and let
$\prod_{k_1,\dots ,k_m}\subset \tilde{\prod}_k$ be the set of all partitions $\sigma\in\tilde{\prod}_k$ such that $|J\cap J'|\leq 1$ for all $J\in\pi $ and all $J'\in\sigma .$ Here $|J|$ is the cardinality of a block $J\in\sigma .$ 
For a partition $\sigma\in\prod_{k_1,\dots ,k_m}$ we define the function $(\otimes_{j=1}^mf_j)_\sigma:B^{|\sigma |}\rightarrow {\mathbb R}$ by replacing all variables of the tensor product $\otimes_{j=1}^mf_j$ that belong to the same block of $\sigma $ by a new common variable, $|\sigma |$ is the number of blocks in $\sigma .$ We denote $\Pi^{(m_{1},\ldots,m_{s}) }_{1,\ldots,s}=\Pi_{1,\ldots,1,\ldots,s,\ldots,s}$, where $i$ repeats $m_{i}$ times for $i=1,\ldots,s$.

For $2 \leq c\leq d$ denote $Y_{c-1}= [0,b]^{d} \times \{ 2b \} \times \{ e_{1},\ldots, e_{c-1} \}$ the space of facets with $d-c+1$ orientations omitted, $\lambda_{c-1}$ is the restriction of measure $\lambda $ onto $Y_{c-1}.$  $\eta_a^{(c-1)}$ is the Poisson process on $Y_{c-1}$ with intensity measure $a\lambda_{c-1}.$ Note that Theorem \ref{lpst} holds also for $\eta_a^{(c-1)}$ and vector $(\tilde{G}_1(\eta_a), \dots ,\tilde{G}_{c-1}(\eta_a)),$ in this case it yields convergence in distribution to $(c-1)$-dimensional Gaussian random vector.
Further for $k\leq p,\,k\leq d,\;k,p\in\N$ let $$Y^p_{[k]}=\{\mathbf x=(u_1,\dots ,u_p)\in Y^p;\; \{j;\,1\leq j\leq d,\;\theta_j( \mathbf{x})>0\}=[k]\}.$$ In the following let $m_j\in \N_0,\; j=1,\dots ,c-1$ and $\sigma\in\Pi^{(m_{1},\ldots,m_{c-1}) }_{1,\ldots,c-1}.$ We will omit arguments and write for short $$\rho_{|\sigma |}=\rho_{|\sigma |}(u_1,\dots ,u_{|\sigma |};\mu_a^{(c)}),$$$$g_\sigma =\left( \otimes_{j=1}^{c-1}\left (\left( g^{(j)} \right)^{\otimes m_{j}}  \right) \right)_{\sigma} (u_{1},\ldots,u_{|\sigma|}).$$ 
\begin{lemma}\label{SimpleLemma}
 There exist $S,R>0$ such that for all $a>0$ \begin{equation}\label{l2}\left\vert\int_{Y^{|\sigma|}}  
g_\sigma\rho_{|\sigma |}\dd\lambda^{|\sigma |}-\int_{Y^{|\sigma|}_{c-1}}  
g_\sigma\dd\lambda_{c-1}^{|\sigma |}\right\vert <Re^{-Sa}.\end{equation}\end{lemma}
\noindent{\bf Proof:} From Lemma \ref{KorelFce} we have $R_1,S>0$ so that $$\sum_{k=1}^{d} \binom{d}{k} \int_{Y^{|\sigma|}_{[k]}}\left( 
\frac{\binom{d-k}{d-c+1}}{\binom{d}{d-c+1}}-R_1e^{-Sa}\right) g_\sigma \dd\lambda^{|\sigma |}<\int_{Y^{|\sigma|}}  
g_\sigma\rho_{|\sigma |}\dd\lambda^{|\sigma |}<$$\begin{equation}\label{ll3}<\sum_{k=1}^{d} \binom{d}{k} \int_{Y^{|\sigma|}_{[k]}}\left( 
\frac{\binom{d-k}{d-c+1}}{\binom{d}{d-c+1}}+R_1e^{-Sa}\right) g_\sigma \dd\lambda^{|\sigma |}.\end{equation}
Since $\binom{n}{k}=0$ for $k>n,$ we have $$\sum_{k=1}^{d} \binom{d}{k}\frac{\binom{d-k}{d-c+1}}{\binom{d}{d-c+1}}=\sum_{k=1}^{c-1}\binom{c-1}{k}$$
and it holds $$\sum_{k=1}^{c-1}\binom{c-1}{k}\int_{Y^{|\sigma|}_{[k]}}g_\sigma\dd\lambda^{|\sigma |}=\int_{Y^{|\sigma|}_{c-1}}g_\sigma\dd\lambda_{c-1}^{|\sigma |}.$$ Finally putting $R=R_1\int_{Y^{|\sigma|}}g_\sigma\dd\lambda^{|\sigma |}$ we obtain the desired bounds for (\ref{l2}) from (\ref{ll3}). \hfill $\Box $

In the following we deal with \begin{equation}\label{mug}\tilde{G}_{j}(\mu^{(c)}_{a}) = \frac{ G_{j}(\mu^{(c)}_{a})- \E G_{j}(\mu^{(c)}_{a})}{a^{j-\frac{1}{2}}},\end{equation} $ 1 \leq j  \leq d$, $2 \leq c \leq d,$ cf. (\ref{stdz}). 
\begin{theorem}\label{CLVT}
For fixed $c=2,\dots ,d$ consider facet processes $\mu_a^{(c)},\; a\geq 1.$ Then
\begin{align}\label{CLV}
(\tilde{G}_{1}(\mu^{(c)}_{a}),\ldots,\tilde{G}_{c-1}(\mu^{(c)}_{a})) \rightarrow M_{c-1},
\end{align}
in distribution as $a\rightarrow\infty ,$ where $M_{c-1}$ is a $(c-1)$-dimensional centered Gaussian random vector with covariance matrix
$\{C_{ij} \}_{i,j=1}^{c-1},\;$ 
$C_{ij}=\langle g_1^{(i)},\,g_1^{(j)}\rangle_1,$ cf. (\ref{ascv}).\end{theorem}
{\bf Proof:} All joint moments of $U$-statistics $G_j(\mu_a^{(c)}),\; j=1,\dots ,c-1,$ can be expressed, see \cite{RefB}, Theorem 2.6, as a finite sum 
\begin{equation}\label{p1}\E \prod_{j=1}^{c-1} G^{m_{j}}_{j}(\mu^{(c)}_{a}) = \sum_{\sigma \in \Pi_{1,\ldots,c-1}^{(m_{1},\ldots,m_{c-1})}}a^{|\sigma|}\int_{Y^{|\sigma|}}
g_\sigma\rho_{|\sigma |}\dd\lambda^{|\sigma |},\end{equation}
while for the Poisson process $\eta_a^{(c-1)}$ all joint moments of $G_j(\eta_a^{(c-1)})$ have form \begin{equation}\label{p2}\E \prod_{j=1}^{c-1} G^{m_{j}}_{j}(\eta^{(c-1)}_{a}) = \sum_{\sigma \in \Pi_{1,\ldots,c-1}^{(m_{1},\ldots,m_{c-1})}}a^{|\sigma|} \int_{Y_{c-1}^{|\sigma|}}g_\sigma\dd\lambda_{c-1}^{|\sigma |},\end{equation} since $\rho_{|\sigma |}(u_1,\dots ,u_{|\sigma |};\eta_a^{(c-1)})\equiv 1.$ 
We deal with the joint moments of the standardized random variables \begin{equation}\label{gjm}\E\prod_{j=1}^{c-1}\tilde{G}_j^{m_j}(\mu_a^{(c)})=\frac{1}{a^q}
\E\prod_{j=1}^{c-1}\left({G}_j(\mu_a^{(c)})-\E {G}_j(\mu_a^{(c)})\right)^{m_j}=\frac{1}{a^q}\sum_{i_1=0}^{m_1}\dots \sum_{i_{c-1}=0}^{m_{c-1}}\end{equation}$$\binom{m_1}{i_1}\dots\binom{m_{c-1}}{i_{c-1}}(-1)^{\sum_{j=1}^{c-1}i_j}\prod_{j=1}^{c-1}\left(\E {G}_j(\mu_a^{(c)})\right)^{i_j}\E\prod_{j=1}^{c-1}{G}_j^{m_j-i_j}(\mu_a^{(c)}),$$ where $q=\sum_{j=1}^{c-1}(j-\frac{1}{2})m_j.$ The terms of the latter formula form products of expressions of type (\ref{p1}), cf. also Theorem 2 in \cite{RefV}, while the terms of the same formula for the Poisson process form products of expressions of type (\ref{p2}). 
 From Lemma \ref{SimpleLemma} we have \begin{equation}\label{gjme}|\E\prod_{j=1}^{c-1}G_j(\mu_a^{(c)})^{m_j-i_j}-\E\prod_{j=1}^{c-1}G_j(\eta_a^{(c-1)})^{m_j-i_j}|\leq C \exp(-Qa)\end{equation} with some constants $C,Q>0.$ All moments involved in these investigations are finite from the assumptions, it suffices to verify this for the Poisson case: $$\E\prod_{j=1}^{c-1}G_j(\eta_a^{(c-1)})^{m_j-i_j}\leq e^{-\lambda (Y)}\sum_{l=1}^\infty\frac{1}{l!}\prod_{j=1}^{c-1}\left[(2b)^{d-j}\binom{l}{j}\right]^{m_j-i_j}a^l\lambda(Y)^l<\infty,$$ using the formula for expectation of a Poisson functional and the upper bound in (\ref{bnds}). The moments of the Poisson $U$-statistics behave asymptotically as the moments of a Gaussian distribution, see \cite{RefLT}, Corollary 4.3 and Proposition 5.1. Because of (\ref{gjme}) the same holds for all joint moments in (\ref{gjm}). Therefore using the moment method we have the central limit theorem (\ref{CLV}) for $(\tilde{G}_{1}(\mu_{a}^{(c)}),\ldots,\tilde{G}_{c-1}(\mu_{a}^{(c)})).$ \hfill $\Box $
\begin{remark}
From Theorem 3 in \cite{RefV}, cf. the beginning of Section \ref{sec}, it follows that the random variables $G_j(\mu_a^{(c)})$ tend to zero in distribution for $j\geq c.$ Therefore it is relevant to study the vector $(\tilde{G}_1(\mu_a^{(c)}),\dots ,\tilde{G}_{c-1}(\mu_a^{(c)}))$ in Theorem \ref{CLVT}.\end{remark}
  
\paragraph{Acknowledgements}
This research was supported by Charles University in Prague, grant SVV 260225, and by the Czech Science Foundation, project 16-03708S.

\end{document}